\input amstex 
\documentstyle{amsppt}
\magnification=1200
\nologo

\redefine\lim{\operatornamewithlimits{lim}}

\define\iso{\tilde\rightarrow}

\define\g{\frak{g}}

\redefine\Im{\operatorname{Im}}
\topmatter  

\title An explicit formula for PBW quantization
\endtitle
\author by Guillermo Corti\~ nas*\endauthor
\affil Departamento de Matem\'atica\\
Facultad de Cs. Exactas y Naturales\\
    Universidad de Buenos Aires\\
       Argentina\endaffil
\thanks (*) Partially supported by grants BID802/OC-AR-PICT 2260 
and UBACyT TW79. This research was carried out while visiting the 
universities of Bielefeld and M\"unster with a Humboldt fellowship.
\endthanks
\address Departamento de
Matem\' atica, 
Ciudad Universitaria Pabell\'on 1,
(1428) Buenos Aires, Argentina.
\endaddress
\email gcorti\@dm.uba.ar\endemail
\leftheadtext{Guillermo Corti\~ nas}
\rightheadtext{An explicit formula for PBW quantization}
\abstract Let $k$ be a field of characteristic zero, $\g$ a $k$-Lie algebra,
$e:S\g@>>>U\g$ the symmetrization map. The PBW quantization is the
one parameter family of associative products:
$$
x\star_t y=\sum_{p=0}^\infty B_p(x,y)t^p\qquad (t\in k)
$$
where $B_p$ is the homogeneous component of degree $-p$ of the map
$B:S\g\otimes_kS\g@>>>S\g$, $B(x,y)=e^{-1}(exey)$. In this paper we give
an explicit formula for $B$. As an application, we prove that for each
$p\ge 0$, $B_p$ is a bidifferential operator of order $\le p$.
\endabstract
\endtopmatter  

\NoBlackBoxes
\document
\head{0. Introduction}\endhead 

We consider (possibly infinite dimensional) Lie algebras over a fixed 
field $k$ of characteristic zero. Let $\g$ be a Lie algebra, $S=S\g$ and 
$U=U\g$ the symmetric and universal enveloping algebras, and 
$\Cal F^0\subset\Cal F^1\subset\dots\subset U$ the coalgebra filtration 
$\Cal F^n:=k+\g+\g^2+\dots+\g^n$. Recall that the Poincar\'e-Birkhoff-Witt 
isomorphism between $S$ and the associated graded ring $G_{\Cal F}U$ is 
induced by the symmetrization map $e:S\iso U$ defined as
$$
e(g_1\dots g_p)=\frac{1}{p!}\sum_{\sigma\in S_p}g_{\sigma(1)}\dots g_{\sigma(p)}
\tag{1}
$$
Thus the associative product
$$
B:S\otimes S@>>>S,\qquad B(x,y)=e^{-1}(exey)\tag{2}
$$
maps $S^{\le n}=\sum_{p=0}^nS^p$ into itself ($n\ge 0$), whence it can be
written as
$$
B=\sum_{p=0}^\infty B_p\tag{3}
$$
where $B_p$ is homogeneous of degree $-p$. We have 
$$
B_0(x,y)=xy,\qquad B_1(x,y)=\frac{1}{2}\{x,y\}\tag{4}
$$
Here $\{,\}$ is the Poisson bracket induced by the Lie bracket of $\g$. 
Thus if $x_1,\dots,x_n$, $y_1,\dots,y_m\in\g$, then
$$
B_1(x,y)=\frac{1}{2}\sum_{i,j}x_1\dots\overset{i}\to{\vee} \dots x_ny_1\dots \overset{j}
\to\vee\dots y_m[x_i,y_j]
$$
In general we will have
$$
B_p(x_1\dots x_n,y_1\dots y_m)=\sum_iz_{i,1}\dots z_{i,m+n-p}\tag{5}
$$
for some elements $z_{ij}\in\g$. In this paper we prove a closed formula
of the form \thetag{5}, and an explicit description of the $z_{i,j}$ as
Lie monomials in the $x_r, y_s$ (Theorem 1.1). As an application of our 
formula, we show that $B_p$ is a bidifferential operator of order $\le p$ 
(Theorem 2.2).

The operators $B_p$ appear naturally in deformation theory and mathematical
physics. One considers
the family of products 
$$
x\star_t y=\sum_{p=0}^\infty B_p(x,y)t^p\qquad (t\in k)
$$
as a one parameter deformation of the usual commutative product $\star_0$ of 
$S$ into the noncommutative product $\star_1$ of the enveloping algebra. If
furthermore $\g$ is finite dimensional, then $S$ can be regarded as the ring
of algebraic functions on the dual $\g^*$, and the $S_t=(S,\star_t)$ as the
rings of functions of a family of noncommutative varieties deforming or
``quantizing'' the Poisson variety $\g^*$. We call this the PBW quantization
because the $B_p$ are defined by means of the Poincar\'e-Birkhoff-Witt theorem.

The idea of the proof of Theorem 1.1 
is to use a well-known expression of the $B_p$ in terms of the 
Campbell-Hausdorff series ([1]) in combination with Dynkin's explicit formula 
for the latter ([3], LA 4.17). Although particular cases of our formula
were known ([1],[2]) this paper is
to our knowledge the first where it appears in its full generality.
An analytic proof of the bidifferentiality of the $B_p$ was given in [1], 
but no estimate of its order is made there. Our proof of the bidifferentiality
is combinatoric and is derived from the explicit formula of Theorem 1.1.

The rest of this paper is organized as follows. The formula for 
$B_p$ is established in section 1 (theorem 1.1). Section
2 is devoted to the proof of its bidifferentiality
(theorem 2.2).
\bigskip
\head{1. A formula for $B_p$}\endhead

In preparation for theorem 1.1 below, we introduce some notation.
Let $n,m\ge 1$, $X=\{x_1,\dots,x_n\}$, $Y=\{y_1,\dots,y_m\}$ two sets of 
noncommuting indeterminates.
If $\alpha\in\{1,\dots,n\}^p$ ($p\ge 1$) is a 
multi-index, then we write:
$$
|\alpha|=p,\quad \Im\alpha=\{\alpha(1),\dots,\alpha(p)\},\qquad ad(x)^\alpha=ad(x_{\alpha(1)})\circ\dots\circ 
ad(x_{\alpha(p)})\tag{6}
$$
The formula of theorem 1.1 below involves the following element of the free
Lie algebra on the disjoint union $X\amalg Y$:
$$
w(X,Y)=\frac{1}{n+m}(w'(X,Y)+w''(X,Y))\tag{7}
$$
Here $w'$ and $w''$ are given by the sums \thetag{8} and \thetag{9};
the restrictions in the summation indexes are explained below; see \thetag{11},
\thetag{12}. 
$$
\align
w'(X,Y)=\sum\frac{(-1)^{p+1}}{p}\frac{ad(x)^{\alpha_1}\circ ad(y)^{\beta_1}\circ\dots\circ ad(x)^{\alpha_p}
(y_k)}{|\alpha_1|!|\beta_1|!\dots |\alpha_p|!}\tag{8}\\
w''(X,Y)=\sum\frac{(-1)^{p+1}}{p}
\frac{ad(x)^{\alpha_1}\circ ad(y)^{\beta_1}\circ\dots\circ ad(y)^{\beta_{p-1}}
(x_k)}{|\alpha_1|!|\beta_1|!\dots |\beta_{p-1}|!}\tag{9}
\endalign
$$
In \thetag{8} the sum is taken over arbitrary $p\ge 1$ and all {\it injective}
multi-indices 
$$
\alpha_i\in\{1,\dots,n\}^{|\alpha_i|}, \qquad 
\beta_j\in\{1,\dots,m\}^{|\beta_j|}\tag{10}
$$
satisfying
$$
\gathered
|\alpha_1|+\dots+|\alpha_p|=n\qquad |\beta_1|+\dots+|\beta_{p-1}|=m-1\\
\quad |\alpha_i|+|\beta_i|\ge 1\quad (1\le i\le p-1)\qquad |\alpha_p|\ge 1\\
\bigcup_{i=1}^p\Im \alpha_i=\{1,\dots,n\}\qquad
\{k\}\cup\bigcup_{i=1}^{p-1}\Im\beta_j=\{1,\dots,m\}
\endgathered\tag{11}
$$
The sum in \thetag{9} is also taken over arbitrary $p\ge 1$, and all injective
multi-indices \thetag{10}, but now
$$
\gathered
|\alpha_1|+\dots+|\alpha_{p-1}|=n-1\qquad |\beta_1|+\dots+|\beta_{p-1}|=m\\
\quad |\alpha_i|+|\beta_i|\ge 1\quad (1\le i\le p-1)\\
\{k\}\cup\bigcup_{i=1}^{p-1}\Im \alpha_i=\{1,\dots,n\}\qquad
\bigcup_{i=1}^{p-1}\Im\beta_j=\{1,\dots,m\}
\endgathered\tag{12}
$$
In the theorem below the we consider the element $w(A,B)$ 
for $A\subset X$, $B\subset Y$. If none of $A$,
$B$ is empty, then $w(A,B)$ is already defined by \thetag{7}; 
we further define
$$
w(\{a\},\emptyset)=w(\emptyset,\{a\})=a\qquad (a\in X\amalg Y)\tag{13}
$$
\bigskip
\definition{Definition 1.0} Let $A_1$, $A_2$ be sets, and $\Cal{P}(A_i)$ the
set of all subsets of $A_i$ ($i=1,2$). A {\it bipartition} of $(A_1,A_2)$ 
is a subset $\pi\subset\Cal{P}(A_1)\times\Cal{P}(A_2)$ such that the following
two conditions are satisfied:
\smallskip
\item{i)} If $S,T\in\pi$ are distinct, then $S_i\cap T_i=\emptyset$ ($i=1,2$).
\smallskip
\item{ii)} $A_i=\cup_{S\in\pi}S_i$ ($i=1,2$)
\smallskip
A bipartition $\pi$ is called {\it special} if $w(S,T)$ is defined for 
all $(S,T)\in\pi$; that is if the following holds
$$
(\emptyset, S) \text{\ \ or\ \ } (S,\emptyset)\in\pi
\Rightarrow \# S=1\tag{14}
$$
\enddefinition
\bigskip
\proclaim{Theorem 1.1} Let $n,m,p\ge 1$,
 $X=\{x_1,\dots,x_n\}$, $Y=\{y_1,\dots,y_m\}$ two sets of indeterminates,
and $B_p$ the operator of \thetag{3} for the free Lie algebra on the 
disjoint union $X\amalg Y$. Then
$$
B_p(x_1\dots x_n, y_1\dots y_m)=\sum_\pi w(\pi_1^X,\pi_1^Y)\dots 
w(\pi_{n+m-p}^X,\pi_{n+m-p}^Y)\tag{15}
$$
Here $w$ is as defined in \thetag{7}, and the sum runs over all special
bipartitions $\pi=\{(\pi_1^X,\pi_1^Y),\dots,(\pi_{n+m-p}^X,\pi_{n+m-p}^Y)\}$
of cardinality $m+n-p$ of $(X,Y)$. In particular,
$$
B_{n+m-1}(x_1\dots x_n,y_1\dots y_m)=w(X,Y)\tag{16}
$$
\endproclaim

\demo{Proof} Write $\g$ for the free Lie algebra on $X\amalg Y$. Let 
$t_1,\dots,t_n$, $u_1,\dots,u_m$ be commuting algebraically
independent variables. Put
$$
x(t)=\sum_{i=1}^nt_ix_i,\qquad y(t)=\sum_{i=1}^mu_iy_i
$$
One checks that $B(x_1\dots x_n,y_1\dots y_m)$ is the coefficient of 
$t_1\dots t_n u_1\dots u_m$ in the product of the exponential series
$$
exp(x(t))exp(y(u))=exp(z(t,u))\tag{17}
$$
Here $z(x(t),y(u))$ is the Campbell-Hausdorff series. Consider the expansion
of $z$ as a series in $t,u$. In order to compute the coefficient of 
$t_1\dots t_n u_1\dots u_m$ in \thetag{17}, all terms in the expansion of
$z$ corresponding to monomials in which any of the $t_i, u_j$ has exponent
$\ge 2$ may be discarded. Each of the remaining terms is an element of $\g$
times a monomial of the form:
$$
t_Au_B:=t_{a_1}\dots t_{{a}_r}u_{{b}_1}\dots u_{{b}_s}
$$
for some subsets $A=\{a_1,\dots,a_r\}\subset\{1,\dots,n\}$, 
$B=\{b_1,\dots, b_s\}\subset\{1,\dots,m\}$. One checks, using Dynkin's formula ([3],LA 4.17),
that the coefficient of $t_Au_B$ in $z$ is precisely the element $w(A,B)$.
The theorem follows from this and the definition of the symmetrization map
\thetag{1}.\qed
\enddemo

In the course of the proof of the theorem above we introduced a notation
which shall be used often in what follows. If $\Cal A$ is a $k$-algebra,
$a_1,\dots,a_n\in\Cal A$ and $S=\{i_1,\dots,i_r\}\subset\{1,\dots,n\}$ is
a subset of $r$ elements, then we write:
$$
a_S:=a_{i_1}\dots a_{i_r}\tag{18}
$$
In particular
$$
a_{\emptyset}=1
$$
\bigskip
\head{2. The bidifferentiality of $B_p$ }\endhead

Let $q,k\ge 0$; define inductively
$$
c_0(q)=1,\qquad c_k(q)=1-\sum_{l=0}^{k-1}c_l(q)\binom{q+k}{k-l}\ \ (k\ge 1)
\tag{19}
$$
\proclaim{Lemma 2.0} Let $\g$ be a Lie algebra, $p\ge 1$, $q\ge 0$, $r\ge 1$, 
$x_1,\dots,x_{p+q}$, $y_1,\dots,y_r\in \g$ and $c_k(q)$ as in \thetag{19}
above. Then, with the notation of \thetag{18}, we have
$$
B_p(x_1\dots x_{p+q},y_1\dots y_r)=\sum_{k=0}^{r-1}c_k(q) (
\sum_{\# S=q+k}
 x_SB_p(x_{{S}^c},y_1\dots y_r))\tag{20}
$$
Here $S\subset\{1,\dots,p+q\}$, and $S^c$ is the complement of $S$. 
The symmetric formula holds for $B_p(y_1\dots y_r,x_1\dots x_{p+q})$.
\endproclaim

\demo{Proof} We may assume that the $x_i,y_j$ are indeterminates and that 
$\g$ is the free Lie algebra. Apply theorem 1.1 to write the left hand side of 
\thetag{20} as a sum of terms indexed by all special bipartitions $\pi$ of 
$(\{1,\dots,p+q\},\{1,\dots,r\})$ of $q+r$ elements. It follows from the 
definition of a special bipartition that the number $d_\pi$ of empty sets
in the list $\pi_1^Y,\dots,\pi_{q+r}^Y$ of subsets of $\{1,\dots, r\}$ is
at least $q$ and at most $q+r-1$. Write $b_k$ for the sum of those terms
whose indexing bipartition has $d_\pi=q+k$ ($0\le k\le r-1$). By definition,
$$
B_p(x_1\dots x_{p+q},y_1\dots y_r)=\sum_{l=0}^{r-1}b_l
$$
Moreover, a counting argument shows that for $0\le k\le r-1$
$$
\sum_{\# S=q+k}x_SB_p(x_{S^c},y_1\dots y_r)=\sum_{i=0}^{r-1-k}\binom{q+k+i}{i}
b_{k+i}
$$
Hence 
$$
\gather
\sum_{k=0}^{r-1}c_k(q)(
\sum_{\# S=q+k}
 x_SB_p(x_{{S}^c},y_1\dots y_r))=
\sum_{k=0}^{r-1}\sum_{i=0}^{r-1-k}c_k(q)
\binom{q+k+i}{i}b_{k+i}\\
=\sum_{l=0}^{r-1}(\sum_{k=0}^l c_k(q)\binom{q+l}{l-k})b_l=
\sum_{l=0}^{r-1}b_l \text{\ \ (by \thetag{19})}\\
=B_p(x_1\dots x_{p+q},y_1\dots y_r)\qed
\endgather
$$
\enddemo

\proclaim{Lemma 2.1 } Let $q\ge 1$, $m\ge 0$. Then:
$$
0= (-1)^m+\sum_{t=1}^q\binom{m+q}{m+t}(-1)^tc_m(t)
$$
\endproclaim
\demo{Proof}Fix $q\ge 1$; the proof is by induction on $m\ge 0$. 
The case $m=0$ is immediate. Assume
$m\ge 1$ and by induction that the lemma holds for $s<m$. Then 
$$\gather
\sum_{t=1}^q\binom{m+q}{m+t}(-1)^tc_m(t)=\\
=\sum_{t=1}^q\binom{m+q}{m+t}(-1)^t-\sum_{t=1}^q\sum_{s=0}^{m-1}(-1)^t
\binom{m+q}{m+t}\binom{t+m}{m-s}c_s(t)\text{\ \ (by \thetag{19})}\\
=(-1)^m\sum_{t=m+1}^{m+q}\binom{m+q}{t}(-1)^t -\sum_{s=0}^{m-1}\sum_{t=1}^q
(-1)^t\binom{m+q}{m-s}\binom{s+q}{s+t}c_s(t)\\
=(-1)^{m+1}(\sum_{t=0}^m(-1)^t\binom{m+q}{t})
+\sum_{s=0}^{m-1}(-1)^{s}\binom{m+q}{m-s}\text{\ \ (by inductive assumption)}\\
=(-1)^{m+1}+(-1)^{m+1}(\sum_{t=1}^m(-1)^t\binom{m+q}{t})+
\sum_{t=1}^m(-1)^{m+t}\binom{m+q}{t}=(-1)^{m+1}\qed
\endgather
$$\enddemo
\medskip
Recall that a $k$-linear endomorphism $F$ of a commutative associative 
algebra $\Cal A$ is a {\it differential operator} of order $\le p$ if
$$
\sum_{S\subset\{1,\dots,p\}}(-1)^S a_{{S}^c}F(a_Sb)=0 \quad 
(a_1,\dots, a_p,b\in\Cal A)\tag{21}
$$
One checks that if $\Cal A$ is generated as a $k$-algebra by a set 
$X\subset\Cal A$
then $F$ satisfies \thetag{21} if and only if it satisfies
$$
\sum_{S\subset\{1,\dots,p+q\}}(-1)^S x_{{S}^c}F(x_S)=0 \quad 
(x_1,\dots, x_{p+q}\in X,\ \ q\ge 0)\tag{22}
$$
\bigskip
\proclaim{Theorem 2.2} Let $\g$ be a Lie algebra, $S=S\g$ the symmetric 
algebra, $a\in S$, $p\ge 1$. Also let $F:S@>>>S$ be one of $B_p(a,\empty)$ or 
$B_p(\empty,a)$. Then $F$ is a differential operator of order $\le p$.
\endproclaim

\demo{Proof}It suffices to show the theorem for $a$ homogeneous. Assume
$a=y_1\dots y_r$, $y_i\in\g$. We shall show that the identity \thetag{22} 
holds for $X=\g$. Put $K=\{1,\dots,p+q\}$. The left hand side of 
\thetag{22} is the sum
of
$$
\sum_{l=p-r+1}^p\sum_{\# S=l}
(-1)^l x_{{K\backslash S}}F(x_S)\tag{23}
$$
and of
$$
\gather
\sum_{t=1}^q\sum_{\# S=p+t}(-1)^{p+t}
x_{K\backslash S}F(x_S)=\\
\sum_{t=1}^q\sum_{\# S=p+t}
\sum_{k=0}^{r-1}\sum_{\gathered T\subset S\\ \# T=p-k\endgathered}
(-1)^{p+t}c_k(t)x_{K\backslash S} x_{S\backslash T}F(x_T)\quad\text{(by 
Lemma 2.0)}\\ 
=\sum_{l=p-r+1}^p\sum_{\# S=l}(\sum_{t=1}^q
\binom{p+q-l}{p+t-l}(-1)^{p+t}c_{p-l}(t))x_{K\backslash S}F(x_S)
\endgather
$$
By lemma 2.1, the sum of the last expression with that of \thetag{23} is zero.
\qed
\enddemo
\subhead{Acknowledgements}\endsubhead I wish to thank Boris Tsygan, from whom
I learned of the relation between the
Campbell-Hausdorff series and the $B_p$. Victor Ginzburg told me that no 
general formula as I was looking for seemed to be known. His comment 
encouraged me to carry out this research, and I am thankful to him for it. 
Thanks also to Alexei Davydov for a useful discussion.


\Refs 

 
\ref\no{1}\by Ginzburg, V.\paper Method of orbits in the representation 
theory of complex Lie groups\jour Funct. Anal. Appl.\vol 15\pages 18-28
\yr 1981
\endref

\ref\no{2}\by Kathotia, V. \paper Konstevich's universal formula for 
deformation quantization and the Campbell-Baker-Hausdorff formula, I
\paperinfo math.QA/9811174 v2 
\endref

\ref\no{3}\by Serre, J.P.\book Lie algebras and Lie groups\publ W. A. Benjamin,
Inc.\yr 1965
\endref


\endRefs
\enddocument